\documentclass[12pt]{amsart}
\DeclareFontFamily{OML}{script}{}
\DeclareFontShape{OML}{script}{m}{it}
{ <5-20> rsfs10 }{}
\DeclareMathAlphabet{\mathscript}{OML}{script}{m}{it}

\renewcommand{\mathcal}[1]{{\mathscript #1}\hspace{0.2ex}}
\usepackage{color}
\ifx\red\undefined
\newcommand{\red}{\color{red}}

\fi
\usepackage[ansinew]{inputenc}
\usepackage[encapsulated]{CJK}

\usepackage{graphicx,graphics}
\usepackage{cite}
\usepackage{pifont}
\usepackage{amsthm}
\usepackage{subfigure}
\usepackage{amscd}
\usepackage{amsmath}
\usepackage{latexsym}
\usepackage{amsfonts}
\usepackage{amssymb}
\usepackage{color}
\usepackage{multicol}
\usepackage{amsmath,amssymb,amsthm,amsfonts,mathrsfs}
\usepackage{hyperref}
\hypersetup{hypertex=true,
            colorlinks=true,
            linkcolor=blue,
            anchorcolor=blue,
            citecolor=blue}
\usepackage{bm}
\allowdisplaybreaks[4]

\newcommand{\Rmnum}[1]{\uppercase\expandafter{\romannumeral #1}}

\ifx\text\undefined
\newcommand{\text}{\mbox}
\fi
\ifx\operatorname\undefined
\newcommand{\operatorname}{\mathop}
\fi

\allowdisplaybreaks[3]

\newcommand\be{\begin{equation}}
\newcommand\ee{\end{euation}}
\newcommand\bea{\begin{eqnarray}}
\newcommand\eea{\end{eqnarray}}
\newcommand\beaa{\begin{eqnarray*}}
\newcommand\eeaa{\end{eqnarray*}}

\newenvironment{eqa}{\begin{equation}%
  \begin{array}{rcl}}{\end{array}\end{equation}}

\newcommand\beqa{\begin{eqa}}
\newcommand\eeqa{\end{eqa}}

\newtheorem{theorem}{Theorem}[section]
\newtheorem{lemma}[theorem]{Lemma}

\newtheorem{remark}[theorem]{Remark}

\theoremstyle{remark}

\newcommand{\void}[1]{}

\def\R{\mathbb{R}}

\numberwithin{equation}{section}

\begin{document}
\begin{CJK}{UTF8}{gkai}

\title[1]{A priori estimates and Liouville-type theorems for the semilinear parabolic equations involving the nonlinear gradient source}
\author{Wenguo Liang and Zhengce Zhang}
\date{\today}
\address[Wenguo Liang]{School of Mathematics and Statistics, Xi'an Jiaotong University,
Xi'an, 710049, P. R. China}
\email{liangwenguo@stu.xjtu.edu.cn}
\address[Zhengce Zhang]{School of Mathematics and Statistics, Xi'an Jiaotong University,
Xi'an, 710049, P. R. China}
\email{zhangzc@mail.xjtu.edu.cn}
\thanks{Corresponding author: Zhengce Zhang}
\thanks{ }

\maketitle

\begin{abstract}
This paper is concerned with the local and global properties of nonnegative solutions for semilinear heat equation $u_t-\Delta u=u^p+M|\nabla u|^q$ in $\Omega\times I\subset \R^N\times \R$, where $M>0$, and $p,q>1$. We first establish the local pointwise gradient estimates when $q$ is subcritical, critical and supercritical with respect to $p$. With these estimates, we can prove the parabolic Liouville-type theorems for time-decreasing ancient solutions. Next, we use Gidas-Spruck type integral methods to prove the Liouville-type theorem for the entire solutions when $q$ is critical. Finally, as an application of the Liouville-type theorem, we use the doubling lemma to derive universal priori estimates for local solutions of parabolic equations with general nonlinearities. Our approach relies on a parabolic differential inequality containing a suitable auxiliary function rather than Keller-Osserman type inequality, which allows us to generalize and extend the partial results of the elliptic equation (Bidaut-V\'{e}ron, Garcia-Huidobro and V\'{e}ron (2020) \cite{veron-sum}) to the parabolic case.\\
\end{abstract}

\section{Introduction}
In this paper, we consider the qualitative properties of parabolic equations with the nonlinear gradient source
\begin{align}\label{eq1}
u_t-\Delta u=u^p+M|\nabla u|^q\ \ {\rm in}\ \Omega\times I,
\end{align}
where $\Omega$ is a domain, $I$ is an interval of $\R$, and $M>0$, $p,q>1$ are parameters. Throughout this paper, the operators $\nabla$ and $\Delta$ only apply to the spatial variables. The solution will always refer to the classical solution contained in $C^{2,1}(\Omega\times I)$.

When $M=0$, equation \eqref{eq1} reduces to
\begin{equation}\label{u-p}
  u_t-\Delta u=u^p.
\end{equation}
Gidas and Spruk used integral estimates in \cite{gidas-spr-1981} to prove the nonexistence of positive global classical stationary solutions when $1<p<p_S:=(N+2)/(N-2)$ for $N>2$. Chen and Li used the moving planes method to provide a new proof of the nonexistence result in \cite{chen-li}. For the steady state equations, the exponent $p_S$ is sharp in the sense that there exist positive bounded radial solutions in $\R^N$ $(N\geq 3)$ when $p\geq p_S$ (see \cite{caff-gidas-spr} and \cite[Section 9]{book-soup}). Then
Bidaut-V\'eron \cite{veron-parab} adapted the integral method to study the blow-up profiles of solutions to parabolic problem \eqref{u-p} for $1<p<p_B$, where
\begin{equation*}
  p_B:=\left\{
  \begin{aligned}
    &\frac{N(N+2)}{(N-1)^2},\quad&\text{if}\ N>1,\\
    &\infty,\quad &\text{if}\ N=1.
  \end{aligned}
  \right.
\end{equation*}
At the same time, Merle and Zaag \cite{Merl-Zaag} proved the nonexistence of positive solutions which satisfy decay assumption of \eqref{u-p} for $1<p<p_S$. However, almost all known results of the optimal Liouville-type theorems for \eqref{u-p} require either a more restrictive assumption on $p$ or deal with a special class of solutions. Recently, Quittner \cite{quitt-optimal} improved the above theorems based on refined energy estimates for the rescaled solution.

Concerning the problems where the nonlinearity contains the spatial derivatives of $u$, the well-known model in this direction is Hamilton-Jacobi equation
\begin{equation}\label{Ham-Job}
  u_t-\Delta u=|\nabla u|^q,
\end{equation}
which arises in both stochastic control problems and KPZ type models of surface growth through ballistic deposition, see  \cite{barles-sto,barles-sto2,kardar-dynamic,krug-surf}.
For the steady state problem, Lions \cite{lions} demonstrated that any classical solution in $\R^N$ for $q>1$ is a constant. The Bernstein estimates are the key to proving this Liouville-type theorem. Moreover, for the half-space steady state problem, Porretta and V\'{e}ron \cite{po-veron} proved the one-dimensional symmetry for $1<q\leq 2$. Filippucci, Pucci, and Souplet obtained a similar symmetry result for $q>2$ in \cite{filip-puc-soup} using the moving planes technique and Bernstein estimates. For the parabolic problem, the Liouville-type theorem for \eqref{Ham-Job} attracted much attention. Souplet and Zhang \cite[Section 3]{soup-zhang} proved gradient estimates for locally upper bounded solution of \eqref{Ham-Job}. Based on these estimates, they demonstrated that ancient solutions remain constants under certain growth conditions at infinity. For further information on the gradient estimates of equation \eqref{Ham-Job}, we refer the reader to \cite{attou,attouchi-sou-20, chang-zhang-23,zhang-hu-10,zhang-li-13}.

For the case of $M<0$ in \eqref{eq1}, the pioneering paper due to Chipot and Weissler \cite{chipot-weis} investigated the possible effects of a gradient term on the global existence or nonexistence, then it was proposed by Souplet \cite{soup-modl} as a population dynamics model to represent the evolution of biological population density. For the case of $M>0$, Tayachi and Zaag \cite{Tayachi} obtained the blow-up behavior of positive solutions. More results concerning the existence and nonexistence, blow-up behavior of solutions to \eqref{eq1} for $M\neq0$ can be found in \cite{veron-chipot,fill-fujita,lu-zhang,zhang-li,zhang19-gradinet}.
It should be noted that few results about gradient estimates of solutions for \eqref{eq1} are presented. However, for the elliptic equation corresponding to \eqref{eq1} when $M>0$, Bidaut-V\'eron, Garcia-Huidobro, and V\'{e}ron \cite{veron-sum} used the Bernstein method and Keller-Osserman's estimates to establish a priori estimates and determine the nonexistence of solutions. By using the integral identity and Young's inequality, Ma, Wu, and Zhang \cite{ma2023liouville} have recently completed the results of nonexistence of solutions in the critical case, that is, $q=2p/(p+1)$.

This paper aims to establish gradient estimates and Liouville-type theorems for equation \eqref{eq1} with  $M>0$. Before proving the main results, we intend to introduce some preliminary notes in this paper.

Let $u$ be a solution of \eqref{eq1}. Define
\begin{equation*}
u_\lambda(x,t)=\lambda^{-1}u(\lambda^{\frac{1-p}{2}}x,\lambda^{1-p}t),\ \ \lambda>0.
\end{equation*}
Clearly, $u_\lambda$ is solution of the equation
\begin{equation*}
u_t-\Delta u=u^p+\lambda^{\frac{(p+1)q-2p}{2}}M|\nabla u|^q.
\end{equation*}
In the case that $\lambda^{(p+1)q/2}\sim\lambda^p$ as $\lambda$ approaches infinity, it holds that $q=2p/(p+1)$. In this context, we refer to $q$ as critical with respect to $p$.
Likewise, if $\lambda^{(p+1)q/2}=o(\lambda^p)$ or $\lambda^{(p+1)q/2}\gg\lambda^p$ as $\lambda$ approaches infinity, we say that $q$ is subcritical or supercritical with respect to $p$, that is, $q<2p/(p+1)$ or $q>2p/(p+1)$. For more details, we refer the reader to \cite{spl-svr}.

As far as we know, there is no Liouville-type theorem for parabolic problem \eqref{eq1} when $M>0$. The presence of the time derivative term and the combined nonlinearities make the problem much more challenging. On the one hand, there has no Keller-Osserman type estimates for parabolic problems, similar to the one for the elliptic problem provided in \cite[Lemma 2.2]{veron-sum}, to study the pointwise gradient estimates. On the other hand, because of the presence of nonlinear gradient term, we cannot find a weighted energy functional for the rescaled problem as in \cite{quitt-optimal}. Thus, the energy estimates technique does not work any more.

Our main goal is to obtain the pointwise and integral gradient estimates of the solutions for \eqref{eq1}. We use the maximum principle for parabolic differential inequality which contains various auxiliary functions to overcome the lack of Keller-Osserman type estimates. In addition, we consider the integral estimates of solutions rather than the energy estimates of the rescaling solutions. It is noteworthy that in local estimates, the gradient term $M|\nabla u|^q$ is regarded as a major or perturbation term. Specifically, we treat the gradient term as the main term when working with the pointwise gradient estimates, and add assumptions on solutions in order to eliminate the nonnegative terms associated with $u^p$. When considering the integral estimates, we regard the terms that contain $M|\nabla u|^q$ as perturbation terms. By applying Young's inequality and attaching suitable assumption on $M$, we can obtain that the perturbation terms will be controlled by other terms.

Now we present the main results. The first result concerns the local pointwise gradient estimates for solutions of \eqref{eq1} in $Q_{T,R}$ when $q$ is subcritical, where
\begin{equation*}
Q_{T,R}=B(x_0,R)\times(0,T),
\end{equation*}
$x_0\in\R^N$, and $R,T>0$.

\begin{theorem}\label{them:q<}
Let $p>1$, $1<q<2p/(p+1)$. For any $M>0$, assume $u$ is the nonnegative solution of \eqref{eq1} in $Q_{T,R}$ that satisfies
\begin{equation}\label{assum-u-t}
u\leq c_{N,p,q}M^{\frac{2}{2p-(p+1)q}}\quad \hbox{and}\quad  u_t\leq 0\quad\text{in}\ Q_{T,R}
\end{equation}
for some constant $c_{N,p,q}>0$.
Then there exists a constant $C:=C(N,M,p,q)>0$, such that
\begin{equation}\label{est-same-nabl}
 |\nabla u(x,t)|\leq  C\left(R^{-1}+R^{-\frac{1}{q-1}}+t^{-\frac{1}{q}}\right)\quad\text{in}\ Q_{T,R/2}.
\end{equation}
\end{theorem}

The gradient estimates of Hamilton-Jacobi equation \eqref{Ham-Job} in \cite[Theorem 3.1]{soup-zhang} is the same as estimates \eqref{est-same-nabl}. It should be noted that if we ignore the coefficient $C(N,M,p,q)$, the power terms in \eqref{est-same-nabl} are independent of $p$. Actually, we first introduce an auxiliary function $f\leq0$ with $f',f''>0$, and $(f''/f')'<0$, and consider the transformation $v=f^{-1}(-u)$. Setting $w=|\nabla v|^2$, we have
\begin{align}\label{int-math-L}
\mathcal L(w):=&w_t-\Delta w-\mathcal H\cdot \nabla w\\\nonumber
\leq&2p(-f)^{p-1}w+2\frac{f''}{(f')^2}(-f)^pw+2\left(\frac{f''}{f'}\right)'w^2-2(q-1)M(f')^{q-2}f''w^{\frac{q+2}{2}}-|D^2v|^2,
\end{align}
where $\mathcal H$ is given in Lemma \ref{lem:Lz}. Note that the last negative term in \eqref{int-math-L} can eliminate the first two nonnegative perturbation terms on the right side of the inequality. Namely, we intend to obtain that
\begin{equation*}\label{ass-u-bound}
2p(-f)^{p-1}w+2\frac{f''}{(f')^2}(-f)^pw-|D^2v|^2\leq0,
\end{equation*}
which requires the assumptions \eqref{assum-u-t} on $u$. It should be pointed out that the boundedness assumption in \eqref{assum-u-t} is the same as the one stated in \cite[Theorem {$\rm A'$}]{veron-sum} for elliptic equations. Therefore, based on the given assumptions, we shall obtain
\begin{equation*}
\mathcal L(w)\leq-2(q-1)M(f')^{q-2}f''w^{\frac{q+2}{2}}.
\end{equation*}
Furthermore, replacing $\mathcal L(w)$ by $\mathcal L(w\eta)$ ensures that the arguments above are always valid for a given cut-off function $\eta$. Thus, it comes that
\begin{equation}\label{math-L-intr}
\mathcal L(w\eta)\leq -\frac{q-1}{2}M(f')^{q-2}f''w^{\frac{q+2}{2}}\eta+C\left(R^{-(q+2)}+R^{-\frac{q+2}{q-1}}\right),
\end{equation}
where $C>0$ depends on $N,M,p$ and $q$. By using the maximum principle to equation \eqref{math-L-intr}, we can conclude that the power terms in the gradient estimates \eqref{est-same-nabl} are independent of $p$.

As a consequence of the local gradient estimates in Theorem \ref{them:q<}, we state the following Liouville-type theorem for bounded solutions in $\R^N\times(-\infty,0)$.

\begin{theorem}\label{them:liou<}
Let $p>1$ and $1<q<2p/(p+1)$. For any $M>0$, assume $u$ is a nonnegative solution of \eqref{eq1} in $\R^N\times(-\infty,0)$ that satisfies
\begin{equation*}
  u\leq c_{N,p,q}M^{\frac{2}{2p-(p+1)q}}\quad \hbox{and}\quad  u_t\leq 0\quad\text{in}\ \R^N\times(-\infty,0)
\end{equation*}
for some constant $c_{N,p,q}>0$.
Then $u\equiv0$.
\end{theorem}

It is remarkable that the time monotonicity condition $u_t\leq 0$ is also required in the proof of the Liouville-type theorem. In fact, we have demonstrated that the ancient solutions are independent of the space variables $x$ based on the estimates \eqref{est-same-nabl}. Consequently, it derives that $0\leq u^p=u_t\leq0$ and thus, $u\equiv0$.

In the case that $q$ is critical, the situation is much more delicate, for the reason that the coefficient of gradient term plays a key role. In what follows, we provide the local gradient estimates for the locally upper bounded solution.
\begin{theorem}\label{them:q=}
Let $p>1$ and $q=2p/(p+1)$. Assume $u$ is a nonnegative solution of \eqref{eq1} in $Q_{T,R}$ that satisfies $u\leq b$ for some constant $b\geq 1$ and
\begin{equation*}
  u_t\leq 0\quad\text{in}\ Q_{T,R}.
\end{equation*}
Then there exists a constant $C:=C(N,M,p)>0$ such that for any
\begin{equation*}
M\geq M_0:=\left(6N(p+1)\right)^{\frac{p}{p+1}}\left(\frac{p+1}{p-1}\right)^{1/2}.
\end{equation*}
the solution satisfies
\begin{equation}\label{three-term}
  |\nabla u(x,t)|\leq Cb\left(R^{-1}+R^{-\frac{p+1}{p-1}}+t^{-\frac{p+1}{2p}}\right)\quad\text{in}\ Q_{T,R/2}.
\end{equation}
\end{theorem}

As a consequence of \eqref{three-term}, the following Liouville-type theorem comes into existence.

\begin{theorem}\label{them:liou=}
Let $p>1$ and $q=2p/(p+1)$. Assume $u$ is a nonnegative solution of \eqref{eq1} in $\R^N\times(-\infty,0)$ for $M\geq M_0$, which satisfies $u_t\leq 0$ in $\R^N\times(-\infty,0)$ and
\begin{equation}\label{u-infty}
  u(x,t)=o\left(|x|+|t|^{\frac{p+1}{2p}}\right)\quad\text{as}\ |x|+|t|^{\frac{p+1}{2p}}\to\infty.
\end{equation}
Then $u\equiv0$. In particular, any bounded nonnegative time-decreasing solution of \eqref{eq1} is $u\equiv0$.
\end{theorem}

Since the gradient estimates in \eqref{three-term} depends on the upper bound of the solution, it is necessary to make a priori assumption on the behavior of $u$ at infinity. However, it is uncertain whether the assumption \eqref{u-infty} is sharp.

Our next results deal with the case that $q$ is supercritical.
\begin{theorem}\label{them:q>}
Let $p>1$, $q>2p/(p+1)$. Assume $u$ is a nonnegative solution of \eqref{eq1} in $Q_{T,R}$, and $u_t\leq\tau$ in $Q_{T,R}$ for any $\tau\geq 0$.

(i) For any $M>0$, the solution $u$ satisfies
\begin{align*}
  |\nabla u(x,t)|\leq C(N,p,q)\left(M^{-\frac{p+1}{(p+1)q-2p}}+M^{-\frac{1}{q}}\tau^{\frac{1}{q}} \right)+C(N,M,p,q)\left(R^{-\frac{1}{q-1}}+t^{-\frac{1}{2(q-1)}}\right)\\\text{in}\  Q_{T,R/2}.
\end{align*}

(ii) For any $M>0$, if
\begin{equation*}
  u\geq c_{N,p,q}\left(M^{-\frac{2}{(p+1)q-2p}}+\tau ^{\frac{1}{p}}\right)\quad \text{in}\ Q_{T,R}
\end{equation*}
for some $c_{N,p,q}>0$. Then  the solution $u$ satisfies
\begin{equation*}
  |\nabla u(x,t)|\leq C\left(R^{-\frac{1}{q-1}}+t^{-\frac{1}{2(q-1)}}\right)\quad\text{in}\ Q_{T,R/2},
\end{equation*}
where $C:=C(N,M,p,q)>0$.
\end{theorem}

\begin{remark}
The upper boundedness assumption for $u_t$ in Theorem \ref{them:q>} can be removed by adding assumptions on $u_0$ and $u$ for the associated Dirichlet and Cauchy problems
\begin{equation*}\label{eq-u0}
  \left\{
\begin{aligned}
  u_t-\Delta u&=u^p+M|\nabla u|^q,\ &&x\in\Omega,\ t>0,\\
  u&=0,&&x\in\partial\Omega,\ t>0,\\
  u(x,t)&=u_0(x),\ &&x\in\Omega.
  \end{aligned}\right.
\end{equation*}
For given $x_0\in \R^N$, we select $R\leq {\rm dist}(x_0,\partial\Omega)$ and then $Q_{T,R}\subset \Omega\times[0,T)$. In fact, assuming that $u$ is uniformly bounded,  $u_0\in C^1(\overline\Omega)$ and $\Delta u_0\in L^\infty(\Omega)$, and using the maximum principle for $u_t$, we have
\begin{equation*}
  u_t\leq \big|\Delta u_0+u_0^p+M|\nabla u_0|^q\big|\quad \text{in}\ \Omega\times(0,T).
\end{equation*}
\end{remark}

We also give a nonexistence result for the ancient solution in what follows.
\begin{theorem}\label{them:liou>}
Let $p>1$ and $q>2p/(p+1)$. Assume $u_t\leq 0$ in $\R^N\times (-\infty,0)$. Then there exists a constant $c_{N,p,q}>0$ such that for any $M>0$, there is no positive solution of \eqref{eq1} in $\R^N\times(-\infty,0)$ that satisfies
\begin{equation}\label{asum-u>}
 u\geq c_{N,p,q}M^{-\frac{2}{(p+1)q-2p}}.
\end{equation}
\end{theorem}

Next, we use the integral Bernstein method to obtain the following result when $q$ is critical and $M$ is small.
\begin{theorem}\label{them:int-q=}
Let $1<p<p_B$ and $q=2p/(p+1)$. Then there exists $\varepsilon_0(N,p)>0$ such that for any $M\in(0,\varepsilon_0)$, the nonnegative solution for \eqref{eq1} in $\R^N\times \R$ is $u\equiv0$.
\end{theorem}

The proof of Theorem \ref{them:int-q=} follows the arguments as in \cite{veron-parab,soupl-dcds-23}. The main difficulty arises from the presence of perturbation term $M|\nabla u|^{q}$, and we have to control the terms $M^2\iint\varphi |\nabla u|^{2q}$, $M\iint\varphi u_t|\nabla u|^q$ and $M\iint\varphi u^p|\nabla u|^q$. Thus, it requires a detailed computation to obtain the integral estimates.

In \cite{sig-07-duke,sig-08-ind}, the authors demonstrated that the Liouville-type theorems ensure the validity of universal estimates. Furthermore, it was found that the Liouville-type theorems were equivalent to the universal estimates in \cite{soupl-dcds-23}. Using the recalling arguments and the doubling property, we will give universal estimates for the solution of \eqref{eq1} as well as the estimates for a more general problem given in \eqref{f-g} below.

\begin{theorem}\label{unifor-esti}
Let $1<p<p_B$, $q=2p/(p+1)$, $T>0$ and $\Omega$ be an arbitrary domain of $\R^N$. Then there exists $\varepsilon_0>0$ such that for any $M\in(0,\varepsilon_0)$, the nonnegative solution of problem \eqref{eq1} in $\Omega\times(0,T)$ satisfies
\begin{equation}\label{uni-2}
  u(x,t)+|\nabla u(x,t)|^{\frac{2}{p+1}}\leq C(N,p)\left(t^{-\frac{1}{p-1}}+(T-t)^{-\frac1{p-1}}+{\rm {dist}}^{-\frac{2}{p-1}}(x,\partial\Omega)\right)\quad\text{in}\ \Omega\times(0,T).
\end{equation}
\end{theorem}

Finally, we will derive universal estimates for general superlinear parabolic equation without scale invariance, given by
\begin{equation}\label{f-g}
  u_t-\Delta u=f(u)+M|\nabla u|^{\frac{2p}{p+1}},
\end{equation}
where $f$ is a nonnegative continuous function satisfying suitable conditions. The result for $q=2p/(p+1)$ in our study is a complement to that obtained in \cite{sig-08-ind,spl-svr} for the case of $u^p+|\nabla u|^q$, where $q$ is considered to be subcritical.

\begin{theorem}\label{them:f-g}
Let $1<p<p_B$, $T>0$ and $\Omega$ be an arbitrary domain of $\R^N$. Assume  $f:[0,\infty)\to \R$ is a continuous function that satisfies
\begin{equation}\label{asm-f}
  \lim_{s\to\infty}s^{-p}f(s)=l>0.
\end{equation}
Then there exists $\varepsilon_1(N,p,l)>0$ such that for any $M\in(0,\varepsilon_1)$, all nonnegative solutions for \eqref{f-g} in $\Omega\times (0,T)$ satisfy
\begin{align}\label{uni}
  u(x,t)+|\nabla u(x,t)|^{\frac{2}{p+1}}\leq C(N,p,f)\left(1+t^{-\frac{1}{p-1}}+(T-t)^{-\frac1{p-1}}+{\rm {dist}}^{-\frac{2}{p-1}}(x,\partial\Omega)\right)\\\nonumber
  \text{in}\ \Omega\times(0,T).
\end{align}
\end{theorem}

Estimates \eqref{uni-2} and \eqref{uni} are universal in the sense that constants $C$ are independent of the solution and domain. These estimates can be written as a more concise form in \cite{sig-08-ind} by using the parabolic distance
\begin{equation*}
  d_P((x,t),(y,s)):=|x-y|+|t-s|^{\frac12}.
\end{equation*}
Namely, denote $D=\Omega\times (0,T)$, $\partial D=(\partial \Omega\times [0,T])\cup \left(\overline\Omega\times\{0,T\}\right)$ the topological boundary of $D$ in $\R^{N+1}$, and $d_P((x,t),\partial D)=\inf_{(y,s)\in \partial D}d_P((x,t),(y,s))$. We will prove the more accurate estimates
\begin{equation}\label{uni-3}
  u(x,t)+|\nabla u(x,t)|^{\frac{2}{p+1}}\leq Cd_P^{-2/(p-1)}((x,t),\partial D)\quad \text{in}\ D
\end{equation}
instead of \eqref{uni-2}, and
\begin{equation}\label{uni-4}
  u(x,t)+|\nabla u(x,t)|^{\frac{2}{p+1}}\leq C\left(1+d_P^{-2/(p-1)}((x,t),\partial D)\right)\quad\text{in}\ D
\end{equation}
instead of \eqref{uni}.

The rest of this paper is organized as follows. In Section \ref{sect-2}, we provide the local pointwise gradient estimates for nonnegative time-decreasing solutions when $p,q>1$, as well as Liouville-type theorems. In Section \ref{sect-3}, we establish integral estimates for the solution and its gradient, along with the Liouville-type theorem for $1<p<p_B$ and $q=2p/(p+1)$. In Section \ref{sect-4}, we use the Liouville-type theorem and doubling method to achieve universal estimates for parabolic equations with a general nonlinearity.

\section{The Bernstein method}\label{sect-2}
In this section, we establish the local pointwise estimates for nonnegative solutions of \eqref{eq1}. These estimates are crucial in proving the Liouville-type theorems for the ancient solution.

Let $\alpha \in(0,1)$ to be chosen later. Set $R^{\prime}=3 R/4$. We select a cut-off function $\eta \in C^2\left(\overline{B}\left(x_0, R^{\prime}\right)\right), 0 \leq \eta \leq 1$, satisfies $\eta=1$ for $|x-x_0|\leq R/2$, $\eta=0$ for $\left|x-x_0\right|=R^{\prime}$ and such that
\begin{equation}\label{dfi-eta}
\left.\begin{array}{rl}
|\nabla \eta| & \leq C R^{-1} \eta^\alpha \\
\left|D^2 \eta\right|+\eta^{-1}|\nabla \eta|^2 & \leq C R^{-2} \eta^\alpha
\end{array}\right\} \text { for }\left|x-x_0\right|<R^{\prime},
\end{equation}
with $C=C(\alpha)>0$. Indeed, such a function is given in \cite{soup-zhang} as $\eta=\rho^k$, where $\rho(x)=1-R'^{-2}|x-x_0|^2$, $x\in\R^N$and $k\geq2/(1-\alpha)$.

We first present a differential inequality of the gradient, which follows the method in \cite{soup-zhang} and will be used several times in this article.
\begin{lemma}\label{lem:Lz}
Let $x_0\in \R^N$, and $R,T>0$. Assume that $u\in C^{2,1}(Q_{T,R})$ is a nonnegative solution of \eqref{eq1}. Define $v=f^{-1}(-u)$, where $f$ is a $C^3$ monotonous function. Set $w=|\nabla v|^2$ and $z=w\eta$. Then at any point such $|\nabla u|>0$, $z$ satisfies
\begin{equation}\label{Lz-0}
\begin{split}
\mathcal L(z)\leq&2p(-f)^{p-1}w\eta+2\frac{f''}{(f')^2}(-f)^{p}w\eta+2\left(\frac{f''}{f'}\right)'w^2\eta\\
&-2(q-1)M(f')^{q-2}f''w^{\frac{q+2}{2}}\eta+\sqrt N|D^2\eta|w+qM(f')^{q-1}|\nabla \eta|w^{\frac{q+1}{2}}\\
&+2\left|\frac{f''}{f'}\right||\nabla \eta|w^{\frac{3}{2}}+C_1(N)|\nabla \eta|^2\eta^{-1}w-|D^2 v|^2\eta\quad\text{in}\  Q_{T,R'},
\end{split}
\end{equation}
where the differential operator $\mathcal L$ is defined by
\begin{equation*}
  \mathcal L(z)=\partial_tz-\Delta z+\mathcal H\cdot\nabla z,
\end{equation*}
and
\begin{equation*}
\mathcal H=\left[qM(f')^{q-1}w^{\frac{q-2}{2}}-2\frac{f''}{f'}\right]\nabla v.
\end{equation*}

\end{lemma}

\noindent {\bf Proof.}
Let $u\in C^{2,1}(Q_{T,R})$ be a nonnegative classical solution of
\begin{equation*}
  u_t-\Delta u=u^p+M|\nabla u|^q.
\end{equation*}
We consider a $C^3$ function $f$ and $f'>0$, the expression of $f$ will be determined later. Throughout the proof, all calculations take place in the region $Q_{T,R'}$. By
\begin{equation}\label{dfi-f}
  u=-f(v),
\end{equation}
we have $u_t=-f'v_t$ and $\nabla u=-f'\nabla v$. Therefore, $v$ satisfies
\begin{equation*}
-f'v_t+f'\Delta v+f''|\nabla v|^2=(-f)^p+M(f')^q|\nabla v|^q.
\end{equation*}
For convenience, the variables of $f, f'$ and $f''$ are omitted here and hereafter. Using the fact that $f'>0$ and setting $w=|\nabla v|^2$, we have the following equation
\begin{equation}\label{eq-v}
  v_t-\Delta v=\frac{(-f)^p}{-f'}-M(f')^{q-1}w^{\frac{q}{2}}+\frac {f''}{f'}w.
\end{equation}
By setting $v_i=\partial v/\partial x_i$ and differentiating \eqref{eq-v} with respect to $x_i$ for $i=1,\ldots,N$, we have
\begin{equation}\label{v_i}
\begin{split}
\partial_tv_{i}-\Delta v_{i}=&p(-f)^{p-1}v_i+\frac{f''}{(f')^2}(-f)^pv_i+\left(\frac{f''}{f'}\right)'v_iw+\frac{f''}{f'}w_i\\
&-M(q-1)(f')^{q-2}f''v_iw^{\frac{q}{2}}-\frac{q}{2}M(f')^{q-1}w^{\frac{q-2}{2}}w_i.
\end{split}
\end{equation}
By the parabolic regularity, we get $\partial_t v_i$, $D^2v_i\in L^r_{loc}(Q_{T,R'})$ for any $r<\infty$. Thus the equations are satisfied in the sense of a.e. in $Q_{T,R'}$. Multiplying \eqref{v_i} by $2v_i$ and summing up with respect to $i$, combined with Bochner's formulate
\begin{equation*}
  2(\nabla v,\nabla(\Delta v))=\Delta w-2|D^2v|^2,
\end{equation*}
we deduce that
\begin{align*}
w_t-\Delta w=&2p(-f)^{p-1}w+2\frac{f''}{(f')^2}(-f)^pw+2\left(\frac{f''}{f'}\right)'w^2+2\frac{f''}{f'}(\nabla w,\nabla v)\\
&-2M(q-1)(f')^{q-2}f''w^{\frac{q+2}{2}}-qM(f')^{q-1}w^{\frac{q-2}{2}}(\nabla w,\nabla v)-2|D^2v|^2.
\end{align*}
We rewrite the equation as
\begin{equation*}
  \mathcal L(w)=-2|D^2v|^2+\mathcal N(w),
\end{equation*}
where the operator $\mathcal L$ is defined by
\begin{equation}\label{math L}
  \mathcal L(w)=\partial_t w-\Delta w+\mathcal H\cdot\nabla w,
\end{equation}

\begin{equation*}
\mathcal N(w):=2p(-f)^{p-1}w+2\frac{f''}{(f')^2}(-f)^pw+2\left(\frac{f''}{f'}\right)'w^2-2(q-1)M(f')^{q-2}f''w^{\frac{q+2}{2}},
\end{equation*}
and
\begin{equation*}
\mathcal H:=\left[qM(f')^{q-1}w^{\frac{q-2}{2}}-2\frac{f''}{f'}\right]\nabla v.
\end{equation*}

Next, we focus on the new function $z=w\eta$. The definition of $\eta$ indicates that $z$ is well-defined in $Q_{T,R'}$. Since
\begin{align}\label{z}\nonumber
  \partial_t(w\eta)=&\eta\partial_tw,\\
  \nabla(w\eta)=&\eta \nabla w+w\nabla \eta,\\\nonumber
  \Delta(w\eta)=&\eta\Delta w+w\Delta\eta+2\nabla\eta\cdot\nabla w,
\end{align}
we can obtain
\begin{align}\label{Lz}\nonumber
\mathcal L(z)=&\partial_t z-\Delta z+\mathcal H\cdot\nabla z\\
=&\eta\mathcal L(w)+w\mathcal L(\eta)-2\nabla w\cdot\nabla \eta\\\nonumber
=&\eta\mathcal N(w)+w\mathcal L(\eta)-2\nabla w\cdot\nabla \eta-2|D^2 v|^2\eta.
\end{align}
Next, we estimate each term in equation \eqref{Lz} on the right side. We start by focus on $w\mathcal L(\eta)$. From Cauchy-Schwartz's inequality, one have
\begin{equation*}
  \frac{f''}{f'}\nabla v\cdot\nabla \eta\leq\left|\frac{f''}{f'}\right|w^{\frac{1}{2}}|\nabla \eta|.
\end{equation*}
By using the fact that for all $\varphi\in C^2(\Omega)$,
 \begin{equation}\label{1/N}
 (\Delta \varphi)^2\leq N|D^2\varphi|^2,
\end{equation}
we obtain
\begin{equation}\label{wL(eta)}
\begin{split}
|w \mathcal L (\eta)|=&w\left|-\Delta \eta+qM(f')^{q-1}w^{\frac{q-2}{2}}\nabla v\cdot\nabla\eta-2\frac{f^{\prime \prime}}{f^{\prime}} \nabla v\cdot \nabla \eta\right|\\
\leq&\sqrt N|D^2\eta|w+qM(f')^{q-1}|\nabla \eta|w^{\frac{q+1}{2}}+2\left|\frac{f''}{f'}\right||\nabla \eta|w^{\frac{3}{2}}.
\end{split}
\end{equation}
From Young's inequality, we know that for any $0<\varepsilon_1(N)<1$, there exists a constant $C(N)>0$ such that
\begin{align}\label{D2v}\nonumber
|\nabla w\cdot\nabla\eta|\leq&|\nabla w||\nabla \eta|\\
\leq&2|D^2v||\nabla v||\nabla\eta|\\\nonumber
\leq& \varepsilon_1|D^2v|^2\eta+C(N)\eta^{-1}|\nabla \eta|^2w.
\end{align}
Substituting \eqref{wL(eta)} and \eqref{D2v} into \eqref{Lz}, and selecting $\varepsilon_1$ suitable small, we obtain
\begin{align*}
\mathcal L(z)\leq&2p(-f)^{p-1}w\eta+2\frac{f''}{(f')^2}(-f)^{p}w\eta+2\left(\frac{f''}{f'}\right)'w^2\eta\\
&-2(q-1)M(f')^{q-2}f''w^{\frac{q+2}{2}}\eta+\sqrt N|D^2\eta|w+qM(f')^{q-1}|\nabla \eta|w^{\frac{q+1}{2}}\\
&+2\left|\frac{f''}{f'}\right||\nabla \eta|w^{\frac{3}{2}}+C(N)|\nabla \eta|^2\eta^{-1}w-|D^2 v|^2\eta,
\end{align*}
which readily implies Lemma \ref{lem:Lz}. \hfill$\Box$

\subsection{The case of $1<q\leq2p/(p+1)$}

In this section, we consider the cases of $1<q<2p/(p+1)$ and $q=2p/(p+1)$, respectively. From the parabolic regularity, we know that $|\nabla u|$ is a locally H\"{o}lder continuous function. Thus $z$ is a continuous function in $\overline{Q}:=\overline{B(x_0,R')}\times[0,T]$. Then, unless $z\equiv0$, there exists a maximum point $(\hat{x},\hat{t})\in \overline{Q}$ such that $z(\hat{x},\hat{t})=\max_{\overline{Q}} z>0$. Since $z=0$ on $\partial B(x_0,R')\times[0,T]$, then $\hat{x}\in B(x_0,R')$. As a result, we have $|\nabla u(x,t)|=f'|\nabla v(x,t)|>0$ in $B(x_0,R')\times(0,T)$.

In the spirit of \cite{attou}, for given $m>0$, we take
\begin{equation}\label{f}
  f(s)=m(s+1)^\gamma-2m,
\end{equation}
where
\begin{equation*}
  \gamma=1+\frac{N}{3(q-1)}.
\end{equation*}
Then $f$ maps $[0,2^{1/\gamma}-1]$ into $[-m,0]$, and $f',f''>0$. We also note that $f'''f'-(f'')^2<0$, thus $\left(f''/f'\right)'<0$.

To apply the maximum principle, we need the following estimates for the right side of \eqref{Lz-0}.

1. Estimates on $\sqrt{N}|D^2\eta|w$. Using Young's inequality, we have
\begin{equation}\label{D^2eta}
  \sqrt{N}|D^2\eta|w\leq \varepsilon_2(q,N)Mm^{q-1}w^{\frac{q+2}{2}}\eta+C(N,M,q)m^{-\frac{2(q-1)}{q}}\left(\eta^{-\frac{2}{q+2}}|D^2\eta|\right)^{\frac{q+2}{q}}.
\end{equation}
The term $|\nabla \eta|^2\eta^{-1}w$ has the similar estimates.

2. Estimates on $qM(f')^{q-1}|\nabla \eta|w^{\frac{q+1}{2}}$. From Young's inequality, we have
\begin{equation}\label{nalba-eta}
  qM(f')^{q-1}|\nabla \eta|w^{\frac{q+1}{2}}\leq \varepsilon_3(q,N)Mm^{q-1}w^{\frac{q+2}{2}}\eta+C(N,M,q)m^{q-1}\left(\eta^{-\frac{q+1}{q+2}}|\nabla \eta|\right)^{q+2}.
\end{equation}

3. Estimates on $\left|f''(f')^{-1}\right||\nabla \eta|w^{3/2}$. Also from Young's inequality, we obtain
\begin{align}\label{w^3/2}\nonumber
  \left|\frac{f''}{f'}\right||\nabla \eta|w^{\frac{3}{2}}\leq& \varepsilon_4(N,M,q)m^{q-1}\left|\frac{f''}{f'}\right|^{\frac{q+2}{3}}w^{\frac{q+2}{2}}\eta+C(N,M,q)m^{-3}\left(\eta^{-\frac{3}{q+2}}|\nabla \eta|\right)^{\frac{q+2}{q-1}}\\
  \leq&\varepsilon_4(N,M,q)C(\gamma,q)m^{q-1}w^{\frac{q+2}{2}}\eta+C(N,M,q)m^{-3}\left(\eta^{-\frac{3}{q+2}}|\nabla \eta|\right)^{\frac{q+2}{q-1}}.
\end{align}

4. Estimates on $|D^2v|^2\eta$. From \eqref{1/N} and \eqref{eq-v}, we derive that
\begin{align}\label{-D2v}
&-|D^2v|^2\\\nonumber
\leq&-\frac{1}{N}(\Delta v)^2\\\nonumber
=&-\frac{1}{N}\left[v_t-\left(\frac{(-f)^p}{-f'}-M(f')^{q-1}w^{\frac{q}{2}}+\frac {f''}{f'}w\right)\right]^2\\\nonumber
=&-\frac{1}{N}\left[v_t^2+\left(\frac{(-f)^p}{-f'}-M(f')^{q-1}w^{\frac{q}{2}}+\frac {f''}{f'}w\right)^2-2v_t\left(\frac{(-f)^p}{-f'}-M(f')^{q-1}w^{\frac{q}{2}}+\frac {f''}{f'}w\right)\right].
\end{align}
Next, following a straightforward calculation, we obtain the expression
\begin{align}\label{b^2}\nonumber
&\left(\frac{(-f)^p}{-f'}-M(f')^{q-1}w^{\frac{q}{2}}+\frac {f''}{f'}w\right)^2\\
=&\frac{(-f)^{2p}}{(f')^2}+M^2(f')^{2q-2}w^q+\left(\frac{f''}{f'}\right)^2w^2+2M(-f)^p(f')^{q-2}w^{\frac{q}{2}}\\\nonumber
&-2(-f)^p\frac{f''}{(f')^2}w-2M(f')^{q-2}f''w^{\frac{q+2}{2}}.
\end{align}
With the assumptions that $u_t\leq 0$ and $f'>0$, we deduce that $v_t=-(f')^{-1}u_t\geq 0$. Hence,
\begin{equation}\label{v_t2}
  -2v_t\left(\frac{(-f)^p}{-f'}-M(f')^{q-1}w^{\frac{q}{2}}\right)=\frac{2v_t}{f'}(u^p+M|\nabla u|^q)\geq 0.
\end{equation}
Using Young's inequality, we obtain that
\begin{equation}\label{v_t1}
  -2v_t\frac{f''}{f'}w\geq -2|v_t|\frac{f''}{f'}w\geq -\frac{1}{4}v_t^2-4\left(\frac{f''}{f'}\right)^2w^2.
\end{equation}
By the assumptions $1<q\leq2p/(p+1)$ and $p>1$, we know that $1<q<2$. Then substituting \eqref{b^2}-\eqref{v_t1} into \eqref{-D2v}, we can obtain
\begin{equation}\label{D^2v-2}
\begin{split}
-|D^2v|^2\leq -&\frac{q-1}{2N}\Bigg[\frac{(-f)^{2p}}{(f')^2}+M^2(f')^{2q-2}w^q+2M(-f)^p(f')^{q-2}w^{\frac{q}{2}}\\
-&2(-f)^p\frac{f''}{(f')^2}w-\left.2M(f')^{q-2}f''w^{\frac{q+2}{2}}-3\left(\frac{f''}{f'}\right)^2w^2\right].
\end{split}
\end{equation}

5. Estimates on $f''(f')^{-2}(-f)^pw\eta$. We consider the second term at the right side of \eqref{Lz-0} and the fourth term at the right side of \eqref{D^2v-2} at the same time. Then it leads to
\begin{align*}
2\left(1+\frac{q-1}{2N}\right)\frac{f''}{(f')^2}(-f)^pw=& 2\left(1+\frac{q-1}{2N}\right)\frac{f''}{f'}\frac{(-f)^p}{f'}w\\
\leq&2\left(1+\frac{q-1}{2N}\right)\frac{\gamma-1}{\gamma}(-f)^{p-1}w\\
\leq&4(-f)^{p-1}w.
\end{align*}

6. Estimates on $(-f)^{p-1}w\eta$. Using Young's inequality and $2(p+2)\leq 6p$, we get
\begin{align}\label{u^pw}\nonumber
6p(-f)^{p-1}w\leq&\frac{1}{q}\left[ \left(\frac{q(q-1)M^2}{2N}\right)^{\frac{1}{q}}(f')^{\frac{2q-2}{q}}w\right]^q\\\nonumber
&+\frac{q-1}{q}\left[ 6p\left(\frac{q(q-1)M^2}{2N}\right)^{-\frac{1}{q}}(-f)^{p-1}(f')^{\frac{2-2q}{q}}\right]^{q/(q-1)}\\
=&(q-1)\frac{M^2}{2N}(f')^{2q-2}w^q+\frac{q-1}{q}(6p)^{\frac{q}{q-1}}\left(\frac{2N}{q(q-1)M^2}\right)^{\frac{1}{q-1}}(-f)^{\frac{(p-1)q}{q-1}}(f')^{-2}\\\nonumber
:=&(q-1)\frac{M^2}{2N}(f')^{2q-2}w^q+C(N,p,q)M^{-\frac{2}{q-1}}(-f)^{\frac{(p-1)q}{q-1}}(f')^{-2}.
\end{align}

Substituting \eqref{D^2eta}-\eqref{w^3/2} and \eqref{D^2v-2}-\eqref{u^pw} into \eqref{Lz-0}, and choosing $\varepsilon_i>0$ suitable small for $i=\{2,3,4\}$,  hence it leads to
\begin{align}\label{Lz-3}
 \mathcal L(z)\leq& -\frac{q-1}{2}M(f')^{q-2}f''w^{\frac{q+2}{2}}\eta+2\left[\left(\frac{f''}{f'}\right)'+\frac{3(q-1)}{N}\left(\frac{f''}{f'}\right)^2\right]w^2\eta\\\nonumber
&+Cm^{-\frac{2(q-1)}{q}}\left(\eta^{-\frac{2}{q+2}}|D^2\eta|\right)^{\frac{q+2}{q}}+Cm^{q-1}\left(\eta^{-\frac{q+1}{q+2}}|\nabla \eta|\right)^{q+2}+Cm^{-3}\left(\eta^{-\frac{3}{q+2}}|\nabla \eta|\right)^{\frac{q+2}{q-1}}\\\nonumber
&+C(N,p,q)M^{-\frac{2}{q-1}}(-f)^{\frac{(p-1)q}{q-1}}(f')^{\frac{2-2q}{q-1}}\eta-\frac{q-1}{2N}\frac{(-f)^{2p}}{(f')^2}\eta\\\nonumber
=&-\frac{q-1}{2}M(f')^{q-2}f''w^{\frac{q+2}{2}}\eta+2\left[\left(\frac{f''}{f'}\right)'+\frac{3(q-1)}{N}\left(\frac{f''}{f'}\right)^2\right]w^2\eta\\\nonumber
&+Cm^{-\frac{2(q-1)}{q}}\left(\eta^{-\frac{2}{q+2}}|D^2\eta|\right)^{\frac{q+2}{q}}+Cm^{q-1}\left(\eta^{-\frac{q+1}{q+2}}|\nabla \eta|\right)^{q+2}+Cm^{-3}\left(\eta^{-\frac{3}{q+2}}|\nabla \eta|\right)^{\frac{q+2}{q-1}}\\\nonumber
&+C(N,p,q)M^{-\frac{2}{q-1}}(f')^{-2}u^{\frac{(p-1)q}{q-1}}\left(1-\frac{q-1}{2N C(N,p,q)}M^{\frac{2}{q-1}}u^{2p-\frac{(p-1)q}{q-1}}\right)\eta.
\end{align}

\noindent\textbf{The proof of Theorem \ref{them:q<}.} Given that $1<q<2p/(p+1)$, it follows that the exponent $2p-(p-1)q/(q-1)<0$. Set
\begin{equation*}
  m=\max_{(x,t)\in \overline{Q}_{T,R}}u(x,t),
\end{equation*}
and select
\begin{equation*}
  m= \left(\frac{2NC(N,p,q)}{q-1}M^{-\frac{2}{q-1}}\right)^{\frac{q-1}{(p+1)q-2p}}:=c_{N,p,q}M^{\frac{2}{2p-(p+1)q}}.
\end{equation*}
It is obvious that the sign of the last term in \eqref{Lz-3} is non-positive. After selecting the suitable function $f$ in \eqref{f}, a calculation shows that
\begin{equation*}
  \left(\frac{f''}{f'}\right)'+\frac{3(q-1)}{N}\left(\frac{f''}{f'}\right)^2=0.
\end{equation*}
From \eqref{Lz-3},  taking $\alpha =3/(q+2)$ in \eqref{dfi-eta} and recalling $\eta\leq 1$, we have
\begin{align}\label{Lz-4}
  \mathcal L(z)\leq& -\frac{q-1}{2}M(f')^{q-2}f''w^{\frac{q+2}{2}}\eta
+C\left(m^{-\frac{2(q-1)}{q}}R^{-\frac{2(q+2)}{q}}+m^{q-1}R^{-(q+2)}+m^{-3}R^{-\frac{q+2}{q-1}}\right)\\\nonumber
\leq&-\frac{\gamma-1}{4}(q-1)Mm^{q-1}w^{\frac{q+2}{2}}\eta+C(N,M,p,q)\left(R^{-(q+2)}+R^{-\frac{q+2}{q-1}}\right),
\end{align}
where we set $C(N,M,p,q)=C\max\{m^{-2(q-1)/q},m^{q-1},m^{-3}\}$ and use the fact that $1<2/q<1/(q-1)$ when $p<2$. Setting
\begin{equation*}
  A_1=C(N,M,m,p,q)\left(R^{-1}+R^{-\frac{1}{q-1}}\right)^2,
\end{equation*}
then we obtain
\begin{equation*}
  \mathcal L(z)\leq -\frac{\gamma-1}{8}(q-1)Mm^{q-1}z^{\frac{q+2}{2}}\quad \text{in}\ \{(x,t)\in Q_{T,R'};\  z(x,t)\geq A_1\}.
\end{equation*}
Next, for a suitable $c=c(N,M,m,p,q)>0$, the function $\psi(t)=ct^{-2/q}$ satisfies
\begin{equation*}
  \psi'(t)\geq -\frac{\gamma-1}{8}(q-1)Mm^{q-1}\psi^{\frac{q+2}{2}}(t).
\end{equation*}
Now we fix $t_0\in (0,T)$ and define $\tilde{z}(x,t)=z(x,t+t_0)-\psi(t)$. It is easy to see that
\begin{equation*}
  \mathcal L(\tilde z)\leq 0\quad\text{in}\ \{(x,t)\in Q_{T-t_0,R'};\ \tilde{z}(x,t)\geq A_1\}.
\end{equation*}
Since $\tilde{z}(x,t)\leq 0$ for sufficiently small $t>0$, then one can apply the maximum principle \cite[Proposition 2.2]{soup-zhang} to deduce that $\tilde{z}(x,t)\leq A_1$, i.e.
\begin{equation*}
  z(x,t+t_0)\leq A_1+\psi(t)\quad\text{in}\ Q_{T-t_0,R'}.
\end{equation*}
Finally, using the fact that $z=w\eta=|\nabla v|^2\eta$ and letting $t_0\to 0$, we obtain
\begin{equation*}
|\nabla v|\leq C(N,M,m,p,q)\left(A_1+t^{-\frac{2}{q}}\right)^{1/2}.
\end{equation*}
Hence by \eqref{dfi-f}, we have
\begin{align}\label{nabla-u}\nonumber
|\nabla u|=&f'|\nabla v|\\
\leq& 2\gamma m C(N,M,m,p,q)\left(A_1+t^{-\frac{2}{q}}\right)^{1/2}\\\nonumber
=& C(N,M,p,q)\left(R^{-1}+R^{-\frac{1}{q-1}}+t^{-\frac{1}{q}}\right)\quad\text{in}\ Q_{T,R/2}.
\end{align}
The proof of Theorem \ref{them:q<} is complete. \hfill$\Box$

\noindent\textbf{The proof of Theorem \ref{them:liou<}.}
Fix $(x_0,t_0)\in \R^N\times(-\infty,0)$. Take $R\geq 1$, $T=R$ and $Q_0=B_R(0)\times(0,T)$. We consider the function $U(y,s):=u(y+x_0,s+t_0-T),\ (y,s)\in Q_0$. From the assumptions on boundedness and monotonicity of the solution with respect to time, we can obtain
\begin{equation*}
  0\leq U\leq c_{N,p,q}M^{\frac{2}{2p-(p+1)q}},
\end{equation*}
and
\begin{equation}
  U_s(y,s)=u_t(y+x_0,s+t_0-T)\leq 0\quad\text{in}\ \overline{Q}_0.
\end{equation}
Applying Theorem \ref{them:q<} to $U$ in $Q_0$, we have
\begin{equation*}
|\nabla u(x_0,t_0)|=|\nabla U(0,T)|\leq C(N,M,p,q)\left(R^{-1}+R^{-\frac{1}{q-1}}+R^{-\frac{1}{q}}\right).
\end{equation*}
It follows that $|\nabla u(x_0,t_0)|=0$ when $R\to\infty$. Since $(x_0,t_0)$ is arbitrary, we obtain that $u$ is independent of space variables $x$. Again using the monotonicity assumption, we get
\begin{equation*}
0\leq u^p=u_t\leq0,
\end{equation*}
thus $u\equiv0$. \hfill$\Box$

\noindent\textbf{The proof of Theorem \ref{them:q=}.}
Suppose $q=2p/(p+1)$, then $2p-(p-1)q/(q-1)=0$. Set $m=b$ and hence $m\geq 1$. From \eqref{u^pw}, we observe that
\begin{equation*}
  C(N,p,q)=\frac{q-1}{q}(6p)^{\frac{q}{q-1}}\left(\frac{2N}{q(q-1)}\right)^{1/(q-1)}.
\end{equation*}
If we select
\begin{equation*}
M\geq M_0:= \left(\frac{2NC(N,p,q)}{q-1}\right)^{\frac{q-1}{2}}=(6N(p+1))^{\frac{p}{p+1}}\left(\frac{p+1}{p-1}\right)^{1/2},
\end{equation*}
then the last term on the right side of \eqref{Lz-3} is non-positive. Using the same arguments as the proof of \eqref{Lz-4} with $m\geq 1$, we get
\begin{equation*}
m^{1-q}\mathcal L(z)\leq -\frac{\gamma-1}{4}(q-1)Mw^{\frac{q+2}{2}}\eta+C(N,M,p,q)\left(R^{-(q+2)}+R^{-\frac{q+2}{q-1}}\right).
\end{equation*}
We set
\begin{equation*}
  A_2=C(N,M,p,q)\left(R^{-1}+R^{-\frac{1}{q-1}}\right)^2
\end{equation*}
and select $\psi(t)=C(N,M,p,q)t^{-2/q}$ satisfying
\begin{equation*}
  \psi'(t)\geq -\frac{\gamma-1}{8}(q-1)M\psi^{\frac{q+2}{2}}(t).
\end{equation*}
Again essentially as before and using the fact that $m^{q-1}\geq1$, we have
\begin{equation*}
  |\nabla v|\leq C(N,M,p)\left(R^{-1}+R^{-\frac{1}{q-1}}+t^{-\frac{1}{q}}\right).
\end{equation*}
Hence, according to $q=2p/(p+1)$ we have
\begin{equation*}
  |\nabla u|=f'|\nabla v|\leq C(N,M,p)b\left(R^{-1}+R^{-\frac{p+1}{p-1}}+t^{-\frac{p+1}{2p}}\right)\quad\text{in}\ Q_{T,R/2}.
\end{equation*}
The proof is thereby complete. \hfill$\Box$

\noindent\textbf{The proof of Theorem \ref{them:liou=}.}
Fix $(x_0,t_0)\in \R^N\times(-\infty,0)$. Take $R\geq 1$ and $T=R^{2p/(p+1)}$. We consider the function $U(y,s)=u(y+x_0,s+t_0-T),\ (y,s)\in Q_0$, where $Q_0=B_R(0)\times(0,T)$. By the assumption \eqref{u-infty}, we know $U\leq b_R$ in $\overline{Q}_0$, where
\begin{equation*}
  b_R:=\sup_{B(x_0,R)\times(t_0-T,t_0)}u=o\left(R+T^{\frac{p+1}{2p}}\right)=o(R)\quad\text{as}\ R\to \infty.
\end{equation*}
Since $U_s\leq0$, then applying Theorem \ref{them:q=} for $U$ and using the fact that $(p+1)/(p-1)>1$, we deduce that
\begin{equation*}
  |\nabla u(x_0,t_0)|=|\nabla U(0,T)|\leq C(N,M,p)R^{-1}b_R.
\end{equation*}
Since $u_t\leq 0$, the result follows from $R$ to $\infty$ and
\begin{equation*}
  0\leq u^p=u_t\leq 0.
\end{equation*}
\hfill$\Box$

\subsection{The case of $q>2p/(p+1)$}
Let $x_0\in \R^N$ be fixed and $R,T>0$. We choose $f(s)=s$, then $v=-u$, $f'=1$ and $f''=0$. Hence, from Lemma \ref{lem:Lz}, we derive
\begin{equation}\label{Lz-I}
  \mathcal L(z)\leq 2pu^{p-1}w\eta+\sqrt{N}|D^2\eta|w+qM|\nabla \eta|w^{\frac{q+1}{2}}+C_1(N)|\nabla \eta|^2\eta^{-1}w-|D^2u|^2\eta\quad\text{in}\ Q_{T,R'}.
\end{equation}
Now we consider each term on the right side of \eqref{Lz-I}. From Young's inequality, we get
\begin{equation*}
  2pu^{p-1}w\eta\leq \frac{1}{4N}u^{2p}\eta+C(N,p)w^{\frac{2p}{p+1}}\eta.
\end{equation*}
Noting that $q> 2p/(p+1)>1$ and using Young's inequality again, we have
\begin{equation*}
  C(N,p)w^{\frac{2p}{p+1}}\eta\leq \frac{M^2}{16N}w^q\eta+C(N,p,q)M^{-\frac{4p}{(p+1)q-2p}}\eta,
\end{equation*}

 \begin{equation*}
  \sqrt{N}|D^2\eta|w\leq \frac{M^2}{16N}w^q\eta+C(N,M,p,q)\left(\eta^{-\frac{1}{q}}|D^2\eta|\right)^{\frac{q}{q-1}},
\end{equation*}

\begin{equation*}
  qM|\nabla \eta|w^{\frac{q+1}{2}}\leq \frac{M^2}{16N}w^{q}\eta+C(N,M,p,q)\left(\eta^{-\frac{q+1}{2q}}|\nabla \eta|\right)^{\frac{2q}{q-1}},
\end{equation*}
and
\begin{equation*}
  C_1(N)|\nabla \eta|^2\eta^{-1}w\leq \frac{M^2}{16N}w^q\eta+C(N,M,p,q)\left(\eta^{-\frac{q+1}{q}}|\nabla \eta|^2\right)^{\frac{q}{q-1}}.
\end{equation*}
Next we consider the term $|D^2u|^2\eta$. From \eqref{1/N}, equation \eqref{eq1} and the assumption $u_t\leq \tau$, we can obtain that
\begin{align*}\label{D2u}\nonumber
-|D^2u|^2\leq& -\frac{1}{N}(\Delta u)^2\\\nonumber
=&-\frac{1}{N}\left[u_t^2+u^{2p}+M^2w^q-2u_t(u^p+Mw^{\frac{q}{2}})+2Mu^pw^{\frac{q}{2}}\right]\\\nonumber
\leq &-\frac{1}{N}u^{2p}-\frac{M^2}{N}w^q+\frac{2\tau}{N}\left(u^p+Mw^{\frac{q}{2}}\right)\\
\leq&-\frac{3}{4N}u^{2p}-\frac{3M^2}{4N}w^q+C(N)\tau^2,
\end{align*}
where we also use Young's inequality for $\tau u^p$ and $\tau w^{q/2}$. Thus, from \eqref{Lz-I} and \eqref{dfi-eta} with $\alpha=(q+1)/(2q)$, we get
\begin{equation}\label{Lz<C1+C2}
  \mathcal L(z)\leq -\frac{M^2}{2N}w^q\eta-\frac{1}{2N}u^{2p}\eta+C(N,p,q)\left(M^{-\frac{4p}{(p+1)q-2p}}+\tau ^2\right)\eta+C(N,M,p,q)R^{-\frac{2q}{q-1}}.
\end{equation}

\noindent\textbf{The proof of Theorem \ref{them:q>}.}
(i) Setting
\begin{equation*}
A_3=C(N,p,q)\left(M^{-\frac{2(p+1)}{(p+1)q-2p}}+M^{-\frac{2}{q}}\tau^{\frac{2}{q}} \right)+C(N,M,p,q)R^{-\frac{2}{q-1}},
\end{equation*}
and then we have
\begin{equation*}
  \mathcal L(z)\leq -\frac{M^2}{4N}z^q\quad\text{in}\ \{(x,t)\in Q_{T,R'};\ \ z(x,t)\geq A_3\}.
\end{equation*}
Next, we select a suitable $c:=c(N,M,q)>0$, such that the function $\psi(t)=ct^{-1/(q-1)}$ satisfies
\begin{equation*}
  \psi'(t)\geq -\frac{M^2}{4N}\psi^q(t).
\end{equation*}
Then again as the proof of Theorem \ref{them:q<}, one gets
\begin{align*}
  |\nabla u|\leq C(N,p,q)\left(M^{-\frac{p+1}{(p+1)q-2p}}+M^{-\frac{1}{q}}\tau^{\frac{1}{q}} \right)+C(N,M,p,q)\left(R^{-\frac{1}{q-1}}+t^{-\frac{1}{2(q-1)}}\right)\\
  \text{in}\ Q_{T,R/2},
\end{align*}
which proves the estimates (i) of Theorem \ref{them:q>}.

(ii) We can select a suitable  constant $c_{N,p,q}>0$, such that
\begin{equation*}
  u^{2p}\geq c^{2p}_{N,p,q} \left(M^{-\frac{2}{(p+1)q-2p}}+\tau ^{\frac{1}{p}}\right)^{2p} \geq C(N,p,q)\left(M^{-\frac{4p}{(p+1)q-2p}}+\tau ^2\right).
\end{equation*}
Then from \eqref{Lz<C1+C2} in $Q_{T,R'}$, we find
\begin{equation*}
  \mathcal L(z)\leq -\frac{M^2}{2N}z^q+C(N,M,p,q)R^{-\frac{2q}{q-1}}.
\end{equation*}
The rest proof is as the precious processes. \hfill$\Box$

\noindent\textbf{The proof of Theorem \ref{them:liou>}.}
When select $\tau =0$, this is an immediate consequence of Theorem \ref{them:q>} (ii). In fact, on the contrary, assume that $u$ is a positive solution for \eqref{eq1} in $\R^N\times(-\infty,0)$, satisfying both \eqref{asum-u>} and $u_t\leq 0$. Then, we can determine that $0<u^p=u_t\leq0$, which is contradictory, essentially as before.  \hfill$\Box$

\section{The integral Bernstein
method}\label{sect-3}
In order to get the integral estimates related with the solution and its gradient, we recall the following lemma from \cite[Lemma 8.9]{book-soup} (see also \cite[Lemma 3.1]{veron-96}). This lemma establishes a set of integral estimates that related to a positive function of its gradient and Laplacian.

\begin{lemma}\label{lem:souplet}
Let $\Omega$ be an arbitrary domain in $\R^N$, $0\leq \varphi\in C_0^\infty(\Omega)$, and $0<v\in C^2(\Omega)$. Fix $a\in \R$, and we denote dy $\int$ the integral on $\Omega$. Set
\begin{equation*}
  I_a=\int \varphi v^{a-2}|\nabla v |^4,\quad J_a=\int\varphi v ^{a-1}|\nabla v |^2\Delta v,\quad K_a=\int \varphi v ^a(\Delta v )^2.
\end{equation*}
Then for any $k\in\R$ with $k\neq-1$, there holds
\begin{equation}\label{souplet-1}
  \alpha I_a+\beta J_a+\gamma K_a\leq \frac{1}{2}\int v ^a|\nabla v |^2\Delta\varphi+\int v ^a\Delta v \nabla v \cdot\nabla \varphi+(a-k)\int v ^{a-1}|\nabla v |^2\nabla v \cdot\nabla \varphi,
\end{equation}
where
\begin{equation*}
  \alpha=-\frac{N-1}{N}k^2+(a-1)k-\frac{a(a-1)}{2},\quad \beta=\frac{N+2}{N}k-\frac{3}{2}a,\quad \gamma=-\frac{N-1}{N}.
\end{equation*}
\end{lemma}

In the following lemma, we shall provide local integral estimates for both the solution and its gradient.

\begin{lemma}\label{est:u^2p}
Let $\Omega$ be an arbitrary domain in $\R^N$, $T>0$, and $0\leq \varphi\in C_0^\infty(\Omega\times (-T,T))$. Let $u$ be a nonnegative solution of \eqref{eq1} in $\Omega\times(-T,T)$. Fix $\theta>0$ and set
\begin{equation*}
  v=u+\theta,\quad f_\theta=u^p-v^p.
\end{equation*}
For given $k\in\R$, denote
\begin{equation*}
I=\iint \varphi v ^{-2}|\nabla v |^4,\quad L=\iint \varphi v ^{2p}, \quad G=M\iint \varphi v ^{-1}|\nabla v |^{2+q},
\end{equation*}
where $\iint$ denotes the integral on $\Omega\times (-T,T)$. Then there holds
\begin{align}\label{I+L-G}\nonumber
\alpha I+\delta L-\beta G\leq& C\iint \varphi\left(M^2|\nabla v |^{2q}+M|v _t||\nabla v |^q+v ^{-1}|v _t||\nabla v |^2+v ^p|v _t|+Mv ^p|\nabla v |^q\right)\\\nonumber
&+C\iint\varphi v_t^2+C\iint |\varphi_t|v ^{p+1}+C\iint |\Delta\varphi||\nabla v |^2\\
&+C\iint|\nabla \varphi|\left(v ^{-1}|\nabla v |^3+\left|v _t-v ^p-M|\nabla v |^q\right||\nabla v |+v ^p|\nabla v |\right)\\\nonumber
&+C\iint \left(\varphi v^{-1}|\nabla v|^2+\varphi v^{p}+|\nabla \varphi||\nabla v|\right)f_\theta+CF_\theta,
\end{align}
where $C=C(N,k)>0$,
\begin{equation*}
  F_\theta=\iint \varphi\left[f_\theta^2-2(v_t-v^p-M|\nabla v|^q)f_\theta\right],
\end{equation*}
and
\begin{align*}
\alpha=&-\frac{N-1}{N}k^2-k,\\
\delta=&-\frac{k}{p}\frac{N+2}{N}-\frac{N-1}{N},\\
\beta=&\frac{N+2}{N}k.
\end{align*}
\end{lemma}

\noindent {\bf Proof.}
We integrate inequality \eqref{souplet-1} with respect to time. Choosing $a=0$ in Lemma \ref{lem:souplet}, we observe that $I=\int_{-T}^{T}I_0$, and set $J=\int_{-T}^{T}J_0$, $K=\int^{T}_{-T}K_0$. It is obtained
\begin{equation*}
\alpha I+\beta J+\gamma K\leq \frac{1}{2}\iint |\nabla v|^2\Delta\varphi+\iint \Delta v\nabla v\cdot\nabla \varphi-k\iint v^{-1}|\nabla v|^2\nabla v\cdot\nabla \varphi,
\end{equation*}
where
\begin{align*}
\alpha =-\frac{N-1}{N}k^2-k,\quad
\beta =\frac{N+2}{N}k,\quad \text{and}\quad \gamma=-\frac{N-1}{N}.
\end{align*}

We next consider the terms $K$ and $J$ by using the perturbation equation
\begin{equation}\label{pertur-eq}
v_t-\Delta v=v^p+M|\nabla v|^q+f_\theta.
\end{equation}
From integrating by parts, we get
\begin{align*}
K=&\iint\varphi(\Delta v)^2\\
=&\iint \varphi\left(v_t^2+v^{2p}+M^2|\nabla v|^{2q}-2v_tv^p-2Mv_t|\nabla v|^q+2Mv^p|\nabla v|^q\right)+F_\theta\\
=&\iint\varphi v_t^2+L+M^2\iint\varphi |\nabla v|^{2q}+\frac{2}{p+1}\iint\varphi_t v^{p+1}\\
&-2M\iint\varphi v_t|\nabla v|^q+2M\iint\varphi v^p|\nabla v|^q+F_\theta,
\end{align*}
and
\begin{align*}
J=&\iint \varphi v^{-1}|\nabla v|^2\left(v_t-v^p-M|\nabla v|^q-f_\theta\right)\\
=&\iint \varphi v^{-1}|\nabla v|^2 v_t-\iint\varphi v^{p-1}|\nabla v|^2-M\iint\varphi v^{-1}|\nabla v|^{2+q}-\iint \varphi v^{-1}|\nabla v|^2f_\theta\\
=&\frac{1}{p}\iint v^p\nabla \varphi\cdot\nabla v+\frac{1}{p}\iint v^p\varphi {\rm{div}} (\nabla v)+\iint\varphi v^{-1}|\nabla v|^2\left(v_t -M|\nabla v|^{q}-f_\theta\right)\\
=&\iint\varphi v^{-1}|\nabla v|^2v_t+\frac{1}{p}\iint v^{p}\nabla v\cdot\nabla \varphi+\frac{1}{p}\iint\varphi v^p v_t-\frac{1}{p}L\\
&-\frac{M}{p}\iint \varphi v^p|\nabla v|^q-M\iint \varphi v^{-1}|\nabla v|^{2+q}-\frac1p\iint\varphi v^p f_\theta -\iint \varphi v^{-1}|\nabla v|^2f_\theta.
\end{align*}
Then substituting the expansions of $K$ and $J$ into \eqref{souplet-1}, we infer that
\begin{align*}
&\alpha I+\delta L-\beta G\\
\leq& C\iint \varphi\left(v_t^2+M^2|\nabla v|^{2q}+M|v_t||\nabla v|^q+v^{-1}|v_t||\nabla v|^2+v^p|v_t|+Mv^p|\nabla v|^q\right)\\
&+C\iint |\varphi_t|v^{p+1}+C\iint |D^2\varphi||\nabla v|^2\\
&+C\iint|\nabla \varphi|\left(v^{-1}|\nabla v|^3+\bigr|v_t-v^p-M|\nabla v|^q\bigr||\nabla v|+v^p|\nabla v|\right)\\
&+C\iint \left(\varphi v^{-1}|\nabla v|^2+\varphi v^{p}+|\nabla \varphi||\nabla v|\right)f_\theta+CF_\theta,
\end{align*}
where
\begin{equation*}
\delta=\gamma-\frac{\beta}{p},
\end{equation*}
which proves the estimates \eqref{I+L-G}. \hfill$\Box$

\begin{remark}
We consider the positive function $u+\theta$ instead of $u$ in order to apply Lemma \ref{lem:souplet}. However, we require the additional complication to handle the perturbation terms $f_\theta$ and $F_\theta$.
\end{remark}

We now derive the integral estimates of solution.

\begin{lemma}\label{u^2p+nabu}
Assume $1<p<p_B$ and $q=2p/(p+1)$. Then there exists $\varepsilon_0>0$ depending on $N$ and $p$, such that for any $M\in (0,\varepsilon_0)$,  the nonnegative solution $u$ of \eqref{eq1} in $Q_R:=B_R\times\left(-R^2,R^2\right)$ satisfies
\begin{equation*}
  \int_{-R^2/2}^{R^2/2}\int_{B_{R/2}}u^{2p}(x,t){\rm d}x{\rm d}t\leq C(N,p)R^{-\frac{4p}{p-1}+N+2}.
\end{equation*}
\end{lemma}

\noindent {\bf Proof.} The proof consists of two steps.

{\em{Step 1: Study of the coefficients and test function.}} Since
\begin{equation*}
  p<p_B=\frac{N(N+2)}{(N-1)^2},
\end{equation*}
we can choose $k<0$ and $k\neq -1$ such that $\alpha,\delta>0$, which is equivalent to
 \begin{equation*}
 \frac{p(N-1)}{N+2} <-k<\frac{N}{N-1}.
 \end{equation*} To estimate the right side of \eqref{I+L-G}, we shall select an appropriate test function $\varphi$. Let $\xi\in C_0^\infty\left(B_1\times(-1,1)\right)$, such that $\xi=1$ in $B_{1/2}\times(-1/2,1/2)$ and $0
\leq\xi\leq 1$. We can fix $\bar\alpha>0$ such that
\begin{equation*}
\frac{3p+1}{4p}<\bar\alpha<1.
\end{equation*}
Fix $R>1$, by taking $\varphi(x,t)=\varphi_R(x,t)=\xi^b(R^{-1}x,R^{-2}t)$ with $b=b(\bar\alpha)>2$ sufficiently large, then we have
\begin{equation}\label{phi-t}
  |\nabla \varphi|\leq CR^{-1}\varphi^{\bar\alpha},\quad |\Delta\varphi|+\varphi^{-1}|\nabla \varphi|^2+|\partial_t\varphi|\leq CR^{-2}\varphi^{\bar\alpha}.
\end{equation}
In the rest of the proof, $\iint$ denotes the integral over $Q_R$.

{\em Step 2: Integral inequalities.}
Now we estimate the right side of \eqref{I+L-G} as follows. Fix $\varepsilon>0$, we first consider the terms containing $\varphi$. By Young's inequality and $q=2p/(p+1)$, we have
\begin{align}\label{M2-u2q}
M^2\varphi |\nabla v|^{2q}\leq& \varepsilon \varphi v^{2p}+C(\varepsilon)M^{\frac4q}\varphi v^{\frac{2(q-2)p}{q}}|\nabla v|^4\\\nonumber
=&\varepsilon \varphi v^{2p}+C(\varepsilon)M^{\frac4q}\varphi v^{-2}|\nabla v|^4.
\end{align}
We also get
\begin{align*}
M\varphi|v_t||\nabla v|^q\leq& \varphi v_t^2+CM^2\varphi |\nabla v|^{2q},\\
\varphi v^{-1}|\nabla v|^2 v_t\leq& \varepsilon \varphi v^{-2}|\nabla v|^4+C(\varepsilon)\varphi v_t^2,\\
\varphi v^p|v_t|\leq&\varepsilon \varphi v^{2p}+C(\varepsilon)\varphi v_t^2,
\end{align*}
and
\begin{equation*}
  M\varphi v^p|\nabla v|^q\leq \varepsilon\varphi v^{2p}+C(\varepsilon)\varphi M^2|\nabla v|^{2q}.
\end{equation*}
Next for the terms containing $|\nabla \varphi|$, we have
\begin{align*}
|\nabla \varphi|v^{-1}|\nabla v|^3\leq & \varepsilon\left(\varphi^{\frac{3}{4}}v^{-\frac{3}{2}}|\nabla v|^3\right)^{\frac{4}{3}}+C(\varepsilon)\left(v^{\frac12}\varphi^{-\frac34}|\nabla \varphi|\right)^4\\
=&\varepsilon \varphi v^{-2}|\nabla v|^4+C(\varepsilon) v^2\varphi^{-3}|\nabla \varphi|^4\\
\leq &\varepsilon \varphi v^{-2}|\nabla v|^4+\varepsilon\left(\varphi^{\frac1p}v^2\right)^p+C(\varepsilon)\left(\varphi^{-\frac{3p+1}{p}}|\nabla\varphi|^4\right)^{\frac{p}{p-1}}\\
=&\varepsilon \varphi v^{-2}|\nabla v|^4+\varepsilon\varphi v^{2p}+C(\varepsilon)\left(\varphi^{-\frac{3p+1}{4p}}|\nabla \varphi|\right)^{\frac{4p}{p-1}},
\end{align*}
\begin{equation}\label{nabla-phi-v_t}
\begin{split}
|\nabla \varphi||v_t||\nabla v|\leq & \varepsilon\left(\varphi^{\frac12}|v_t|\right)^2+C(\varepsilon)\left(\varphi^{-\frac12}|\nabla \varphi||\nabla v|\right)^2\\
  =&\varepsilon\varphi v_t^2+C(\varepsilon)\varphi^{-1}|\nabla \varphi|^2|\nabla v|^2\\
  \leq& \varepsilon\varphi v_t^2+\varepsilon\left(v^{-1}|\nabla v|^2\varphi^{\frac12}\right)^2+C(\varepsilon)\left(v|\nabla \varphi|^2\varphi^{-\frac32}\right)^2\\
  =&\varepsilon\varphi v_t^2+\varepsilon \varphi v^{-2}|\nabla v|^4+C(\varepsilon)v^2\varphi^{-3}|\nabla \varphi|^4\\
  \leq &\varepsilon\varphi v_t^2+\varepsilon \varphi v^{-2}|\nabla v|^4+\varepsilon\varphi v^{2p}+C(\varepsilon)\left(\varphi^{-\frac{3p+1}{4p}}|\nabla \varphi|\right)^{\frac{4p}{p-1}},
\end{split}
\end{equation}
and
\begin{align*}
|\nabla \varphi|v^p|\nabla v|\leq &\varepsilon\left(\varphi^{\frac14}v^{-\frac12}|\nabla v|\right)^4+C(\varepsilon)\left(\varphi^{-\frac14}v^{p+\frac12}|\nabla\varphi|\right)^{\frac43}\\
=&\varepsilon\varphi v^{-2}|\nabla v|^4+C(\varepsilon)\varphi^{-\frac13}|\nabla \varphi|^{\frac43}v^{\frac{4p+2}{3}}\\
\leq &\varepsilon\varphi v^{-2}|\nabla v|^4+\varepsilon\left(\varphi^{\frac{2p+1}{3p}}v^{\frac{4p+2}{3}}\right)^{\frac{3p}{2p+1}}+C(\varepsilon)\left(\varphi^{-\frac{3p+1}{3p}}|\nabla \varphi|^{\frac43}\right)^{\frac{3p}{p-1}}\\
=&\varepsilon\varphi v^{-2}|\nabla v|^4+\varepsilon \varphi v^{2p}+C(\varepsilon)\left(\varphi^{-\frac{3p+1}{4p}}|\nabla\varphi|\right)^{\frac{4p}{p-1}}.
\end{align*}
Since $q=2p/(p+1)$, we obtain
\begin{align*}
M|\nabla \varphi||\nabla v|^{q+1}\leq& CM|\nabla \varphi||\nabla v|\left[\left(v^{-\frac{q}{2}}|\nabla v|^q\right)^{\frac2q}+v^{\frac{q}{2-q}}\right]\\
=&CM|\nabla \varphi|v^{-1}|\nabla v|^3+CM|\nabla \varphi|v^p|\nabla v|.
\end{align*}
We now consider the terms containing $|\varphi_t|$ and $|\Delta\varphi|$. By Young's inequality, we obtain
\begin{equation}\label{phi-u-p+1}
\begin{split}
|\varphi_t|v^{p+1}\leq& \varepsilon\left(\varphi^{\frac{p+1}{2p}}v^{p+1}\right)^{\frac{2p}{p+1}}+C(\varepsilon)\left(\varphi^{-\frac{p+1}{2p}}|\varphi_t|\right)^{\frac{2p}{p-1}}\\
=&\varepsilon\varphi v^{2p}+C(\varepsilon)\left(\varphi^{-\frac{p+1}{2p}}|\varphi_t|\right)^{\frac{2p}{p-1}}
\end{split}
\end{equation}
and
\begin{align*}
|\Delta\varphi||\nabla v|^2\leq &\varepsilon \left(\varphi^{\frac12}v^{-1}|\nabla v|^2\right)^2+C(\varepsilon)\left(\varphi^{-\frac12}v|\Delta\varphi|\right)^2\\
=&\varepsilon \varphi v^{-2}|\nabla v|^4+C(\varepsilon) v^2\varphi^{-1}|\Delta\varphi|^2\\
\leq &\varepsilon \varphi v^{-2}|\nabla v|^4+\varepsilon \left(v^2\varphi^{\frac1p}\right)^p+C(\varepsilon)\left(\varphi^{-\frac{p+1}{p}}|\Delta\varphi|^2\right)^{\frac{p}{p-1}}\\
=&\varepsilon \varphi v^{-2}|\nabla v|^4+\varepsilon \varphi v^{2p}+C(\varepsilon)\left(\varphi^{-\frac{p+1}{2p}}|\Delta\varphi|\right)^{\frac{2p}{p-1}}.
\end{align*}
It remains to consider the terms containing $f_\theta$ and $F_\theta$. By Young's inequality, we have
\begin{equation*}
  \iint\varphi v^{-1}|\nabla v|^2f_\theta\leq \varepsilon \iint\varphi v^{-2}|\nabla v|^4+C(\varepsilon)\iint \varphi f^2_\theta,
\end{equation*}
\begin{equation*}
  \iint \varphi v^p f_\theta\leq \varepsilon \iint\varphi v^{2p}+C(\varepsilon)\iint\varphi f^2_\theta,
\end{equation*}
and
\begin{align*}
&\iint|\nabla \varphi||\nabla v|f_\theta\\
&\leq \varepsilon \iint\left(\varphi^{\frac14}v^{-\frac12}|\nabla v|\right)^4+C(\varepsilon)\iint\left(v^{\frac12}\varphi^{-\frac14}|\nabla \varphi|f_\theta\right)^{\frac43}\\
&\leq \varepsilon\iint\varphi v^{-2}|\nabla v|^4+\varepsilon \iint\left(\varphi^{\frac1{3p}} v^{\frac23}\right)^{3p}+C(\varepsilon)\iint\left(\varphi^{-\frac{p+1}{3p}}|\nabla \varphi|^{\frac43}f_\theta^{4/3}\right)^{\frac{3p}{3p-1}}\\
&=\varepsilon\iint\varphi\left( v^{-2}|\nabla v|^4+ v^{2p}\right)+C(\varepsilon)\iint\varphi^{-\frac{p+1}{3p-1}}|\nabla \varphi|^{\frac{4p}{3p-1}}f_\theta^{4p/(3p-1)}\\
&\leq \varepsilon\iint \varphi\left( v^{-2}|\nabla v|^4+ v^{2p}\right)+C(\varepsilon)\iint\left(\varphi^{-(3p+1)}|\nabla \varphi|^{4p}\right)^{\frac1{p-1}}+C(\varepsilon)\iint \varphi f_\theta^2.
\end{align*}
According to the definition of $F_\theta$, we find
\begin{equation*}
F_\theta\leq\iint\varphi\left(v_t^2+\varepsilon v^{2p}+\varepsilon M^2|\nabla v|^{2q}+C(\varepsilon)f^2_\theta\right).
\end{equation*}
Therefore, substituting these estimates into \eqref{I+L-G}, we derive
\begin{equation}\label{I+L<R}
\begin{split}
&\alpha I+\delta L-\beta G\\
\leq
&C\left(\varepsilon I+\varepsilon L+C(\varepsilon)M^{\frac4q}I +\iint\varphi v_t^2\right)+C(\varepsilon)\iint\varphi f^2_\theta\\ &+C(\varepsilon)\left[\iint\left(\varphi^{-\frac{3p+1}{4p}}|\nabla \varphi|\right)^{\frac{4p}{p-1}}
+\iint\left(\varphi^{-\frac{p+1}{2p}}|\Delta\varphi|\right)^{\frac{2p}{p-1}}+\iint\left(\varphi^{-\frac{p+1}{2p}}|\varphi_t|\right)^{\frac{2p}{p-1}}\right].
\end{split}
\end{equation}
Finally, we estimate the time derivative term of $v$. From \eqref{pertur-eq}, we get
\begin{align}\label{iint v_t}
&\iint \varphi v_t^2\\\nonumber
=&\iint\varphi v_t\left(v^p+M|\nabla v|^q+\Delta v+f_\theta\right)\\\nonumber
  =&\frac{1}{p+1}\iint \varphi\left(v^{p+1}\right)_t+M\iint\varphi v_t|\nabla v|^q+\iint \varphi v_t\Delta v+\iint\varphi v_tf_\theta\\\nonumber
  \leq& \frac1{p+1}\iint |\varphi_t|v^{p+1}+\frac{1}2\iint\left(\varphi v_t^2+M^2\varphi|\nabla v|^{2q}\right)-\iint \left(\nabla \varphi v_t+\varphi\nabla v_t \right)\cdot\nabla v+\iint\varphi v_tf_\theta\\\nonumber
  \leq &C(p)\iint |\varphi_t|v^{p+1}+\frac{1+2\varepsilon}2\iint\left(\varphi v_t^2+M^2\varphi|\nabla v|^{2q}\right)+\iint\varphi\partial_t\left(\frac{|\nabla v|^2}{2}\right)\\\nonumber
  &+\varepsilon\iint \varphi v^{-2}|\nabla v|^4+\varepsilon\iint\varphi v^{2p}+C(\varepsilon)\iint\left(\varphi^{-\frac{3p+1}{4p}}|\nabla \varphi|\right)^{\frac{4p}{p-1}}+C(\varepsilon)\iint\varphi f^2_\theta.
\end{align}
Integrating by parts and using Young's inequality, we deduce that
\begin{equation}\label{part-t-u}
\begin{split}
\iint \varphi\partial_t\left(\frac{|\nabla v|^2}{2}\right)=&-\frac12\iint\varphi_t|\nabla v|^2\\
 \leq&\varepsilon\iint\varphi v^{-2}|\nabla v|^4+\varepsilon\iint\varphi v^{2p}+C(\varepsilon)\iint\left(\varphi^{-\frac{p+1}{2p}}|\varphi_t|\right)^{\frac{2p}{p-1}}.
\end{split}
\end{equation}
Substituting \eqref{M2-u2q}, \eqref{phi-u-p+1} and \eqref{part-t-u} into \eqref{iint v_t}, we have that
\begin{equation}\label{int-phi-t}
\begin{split}
  \iint \varphi v_t^2\leq& C\left(\varepsilon I+\varepsilon L+ C(\varepsilon)M^{\frac4q}I\right)\\
  &+C(\varepsilon)\iint\left[\left(\varphi^{-\frac{3p+1}{4p}}|\nabla \varphi|\right)^{\frac{4p}{p-1}}+\left(\varphi^{-\frac{p+1}{2p}}|\varphi_t|\right)^{\frac{2p}{p-1}}\right]+C(\varepsilon)\iint\varphi f^2_\theta.
\end{split}
\end{equation}
Hence, from \eqref{int-phi-t} and \eqref{I+L<R}, we have
\begin{equation}\label{est-I-L-G}
\begin{split}
&\alpha I+\delta L-\beta G\\
 \leq& C\left(\varepsilon I+\varepsilon L+C(\varepsilon)M^{\frac4q}I\right)+C(\varepsilon)\iint\varphi f^2_\theta\\
  &+C(\varepsilon)\iint\left[\left(\varphi^{-\frac{3p+1}{4p}}|\nabla \varphi|\right)^{\frac{4p}{p-1}}+\left(\varphi^{-\frac{p+1}{2p}}|\Delta\varphi|\right)^{\frac{2p}{p-1}}+\left(\varphi^{-\frac{p+1}{2p}}|\varphi_t|\right)^{\frac{2p}{p-1}}\right].
\end{split}
\end{equation}
We select $\varepsilon=\varepsilon(N,p)>0$ suitable small, and there exists $\varepsilon_0(N,p)>0$ such that, for $0<M\leq \varepsilon_0$,
\begin{equation*}
  C(\varepsilon)M^{\frac4q}=C(N,p)M^{\frac4q}\leq \frac{\alpha}{2}.
\end{equation*}
Combining with \eqref{phi-t}, and using the fact $\bar\alpha>(3p+1)/(4p)$, we obtain
\begin{equation*}
    \frac{\alpha}{2}I+\frac{\delta}{2} L-\beta G\leq C(N,p)R^{-\frac{4p}{p-1}+N+2}+C(\varepsilon)\iint\varphi f^2_\theta.
\end{equation*}
Using the fact that $\beta<0$ and
\begin{equation*}
  \int_{-R^2/2}^{R^2/2}\int_{B_{R/2}}\varphi v^{2p}(x,t){\rm d}x{\rm d}t\to \int_{-R^2/2}^{R^2/2}\int_{B_{R/2}} u^{2p}(x,t){\rm d}x{\rm d}t\quad \text{as}\ \theta\to0,
\end{equation*}
then Lemma \ref{u^2p+nabu} follows as $\theta\to 0$. \hfill$\Box$

\noindent\textbf{The proof of Theorem \ref{them:int-q=}.}
From Lemma \ref{u^2p+nabu}, we know that
\begin{equation*}
  \int_{-R^2/2}^{R^2/2}\int_{B_{R/2}}u^{2p}(x,t){\rm d}x{\rm d}t\leq C(N,p)R^{-\frac{4p}{p-1}+N+2}.
\end{equation*}
Since $p<p_B\leq p_S$, letting $R\to \infty$, we deduce $u\equiv0$. \hfill$\Box$

\begin{remark}
Since equation \eqref{eq1} is scaling invariant when $q=2p/(p+1)$. Consequently, Theorem \ref{them:int-q=} can be proved using the universal bound of the solution as in \cite[Theorem 21.2]{book-soup}. In fact, by repeating the same argument as in Lemma \ref{u^2p+nabu} (replacing the integration region $Q_R$ by $Q_1$), we obtain that
\begin{equation}\label{int-bound}
  \int_{-1/2}^{1/2}\int_{B_{1/2}} u^{2p}(x,t){\rm d}x{\rm d}t\leq C(N,p).
\end{equation}
For any $R>0$, then the rescaling function
\begin{equation*}
u_R(x,t)=R^{\frac2{p-1}}u(Rx,R^2t)
\end{equation*}
is also a solution of equation \eqref{eq1} in $\R^N\times \R$ for $q=2p/(p+1)$. Then from \eqref{int-bound}, we have
\begin{align*}
 \int_{-R^2/2}^{R^2/2}\int_{B_{R/2}} u^{2p}(y,s){\rm d}y{\rm d}s&=R^{-\frac{4p}{p-1}+N+2}\int_{-1/2}^{1/2}\int_{B_{1/2}} u_R^{2p}(x,t){\rm d}x{\rm d}t\\
 &\leq C(N,p)R^{-\frac{4p}{p-1}+N+2}.
\end{align*}
Letting $R\to\infty$, we deduce $u\equiv0$.
\end{remark}

\section{Universal estimates via Liouville-type theorems}\label{sect-4}

We first recall the useful doubling lemma due to Pol\'{a}\v{c}ik, Quitter and Souplet \cite[Lemma 5.1]{sig-07-duke}, which allows us to structure the rescaling procedure and obtain the local estimates. The doubling property is an expansion derived from Hu \cite{Hu-96-doubing}.

\begin{lemma}\label{lem:doub}
Let $(\mathscr{X}, d)$ be a complete metric space with metric $d$, and let $\emptyset \neq D \subset \Sigma \subset \mathscr{X}$, with $\Sigma$ closed. Set $\Gamma=\Sigma \backslash D$. Let $M: D \rightarrow(0, \infty)$ be bounded on compact subsets of $D$ and fix a real number $k>0$. If $y \in D$ is such that
\begin{equation*}
  M(y){\rm{dist}}(y,\Gamma)>2k,
\end{equation*}
then there exists $x\in D$ such that
\begin{equation*}
  M(x){\rm{dist}}(x,\Gamma)>2k,\ \ M(x)\geq M(y),
\end{equation*}
and
\begin{equation*}
  M(z)\leq 2M(x)\ \ \text{\rm for all}\ z\in D\cap\overline{B}_{\mathscr{X}}\left(x,kM^{-1}(x)\right).
\end{equation*}
\end{lemma}

\noindent\textbf{The proofs of Theorems \ref{unifor-esti} and \ref{them:f-g}.}  We set
\begin{equation*}
  M(u)=u^{\frac{p-1}{2}}+|\nabla u|^{\frac{p-1}{p+1}},
\end{equation*}
then the estimates \eqref{uni-3} and \eqref{uni-4} can be written as
\begin{equation}\label{est-sigma}
  M(u(x,t))\leq C\left(\sigma+d_P^{-1}((x,t),\partial D)\right),\quad (x,t)\in D,
\end{equation}
where $C>0$ and
\begin{equation*}
  \sigma:=\left\{
\begin{aligned}
  &0, &&\text{under the assumptions of Theorem \ref{unifor-esti}},\\
  &1, &&\text{under the assumptions of Theorem \ref{them:f-g}}.
\end{aligned}
\right.
\end{equation*}
Assume that the estimates \eqref{est-sigma} fails. Then there exist sequences of domains $\Omega_k$, times $T_k>0$, solutions $u_k$ of \eqref{eq1} or \eqref{f-g} in $D_k:=\Omega_k\times (0,T_k)$, respectively, and points $(y_k,\tau_k)\in D_k$ such that
\begin{equation}\label{M>2k}
M_k(y_k,\tau_k):=M(u_k(y_k,\tau_k))>2k\left(\sigma+d_P^{-1}\left((y_k,\tau_k),\partial D_k\right)\right)\geq 2kd_P^{-1}((y_k,\tau_k),\partial D_k).
\end{equation}
We use Lemma \ref{lem:doub} with $\mathscr{X}=\R^{N+1}$, equipped with the parabolic distance $d_P$, $\Sigma=\Sigma_k=\overline{D}_k$ and $\Gamma=\partial D_k$. Then there exist points $(x_k,t_k)\in D_k$, such that
\begin{equation}\label{M-xk}
  M_k(x_k,t_k)>2kd_P^{-1}((x_k,t_k),\partial D_k),\quad M_k(x_k,t_k)\geq M_k(y_k,\tau_k)>2k\sigma,
\end{equation}
and
\begin{equation}\label{M<2Mk}
  M_k(x,t)\leq 2M_k(x_k,t_k)\quad \text{in}\ \hat D_k:=\left\{(x,t)\in\R^{N+1};\ d_P\left((x,t),(x_k,t_k)\right)\leq kM_k^{-1}(x_k,t_k)\right\}.
\end{equation}
Using the first inequality of \eqref{M-xk}, we have
\begin{equation*}
  d_P((x,t),(x_k,t_k))<\frac12 d_P((x_k,t_k),\partial D_k),\quad(x,t)\in \hat D_k,
\end{equation*}
which implies $\hat D_k\subset D_k$. We set
\begin{equation}\label{lamb}
  \lambda_k=M_k^{-1}(x_k,t_k),
\end{equation}
and then rescale $u_k$ by setting
\begin{equation}\label{scal}
  v_k(y,s)=\lambda _k^{2/(p-1)}u_k\left(x_k+\lambda_ky,t_k+\lambda_k^2s\right),\quad (y,s)\in \tilde D_k,
\end{equation}
where
\begin{equation*}
  \tilde D_k:=\left\{y\in \R^N;\ |y|< k/2\right\}\times \left(-k^2/4,k^2/4\right).
\end{equation*}
Then the function $v_k$ satisfies
\begin{equation*}
  \partial_sv_k-\Delta v_k=v^p_k+M|\nabla v_k|^{2p/(p+1)}\quad\text{in}\ \tilde D_k,
\end{equation*}
or
\begin{equation*}
 \partial_sv_k-\Delta v_k=\lambda_k^{2p/(p-1)}f\left(\lambda_k^{-2/(p-1)}v_k\right)+M|\nabla v_k|^{2p/(p+1)}\quad \text{in}\ \tilde D_k
\end{equation*}
under the assumptions of Theorem \ref{unifor-esti} or Theorem \ref{them:f-g}, respectively. From \eqref{lamb} and \eqref{scal}, we have
\begin{equation*}
  \left[v_k^{(p-1)/2}+|\nabla v_k|^{(p-1)/(p+1)}\right](0,0)=\lambda_k M_k(x_k,t_k)=1.
\end{equation*}
By \eqref{M<2Mk}, we also have
\begin{equation*}
\left[v_k^{(p-1)/2}+|\nabla v_k|^{(p-1)/(p+1)}\right](y,s)\leq 2,\quad (y,s)\in \tilde D_k.
\end{equation*}
Using the assumption \eqref{asm-f} on $f$ and the boundedness of $\nabla v_k$, we derive
\begin{equation*}
0\leq\lambda_k^{2p/(p-1)}f\left(\lambda_k^{-2/(p-1)}v_k(y,s)\right)+M|\nabla v_k(y,s)|^{2p/(p+1)}\leq C(1+M),\quad (y,s)\in \tilde D_k.
\end{equation*}
where $C=C(N,p,l)>0$.

Moreover, by the parabolic $L^p$ estimates and embedding theorem,  up to a subsequence, there exists a $v\in W^{2,1;r}_{\rm loc}(\R^N\times\R)$ such that
\begin{equation*}
  v_k\to v\quad \text{in}\  C^{2,1}_{\rm loc}(\R^N\times\R)
\end{equation*}
for any $1<r<\infty$ and $v$ is nontrivial. Thus, under the assumptions of Theorem \ref{unifor-esti}, $v$ is a nonnegative nontrivial solution of equation \eqref{eq1} in $\R^N\times \R$ with $0<M\leq \varepsilon_0$, this contradicts to Theorem \ref{them:int-q=} and thus proves Theorem \ref{unifor-esti}.

Now we consider the case of Theorem \ref{them:f-g}. Applying the maximum principle, there exists a $s^*\in [-\infty,0)$ such that
\begin{equation*}
v=0\quad\text{\rm in}\ \R^N\times (-\infty,s^*],
\end{equation*}
and
\begin{equation*}
v>0\quad\text{\rm in}\ \R^N\times (s^*,\infty).
\end{equation*}
Using \eqref{M-xk} and \eqref{lamb}, we obtain
\begin{equation*}
  \lambda_k\to 0\quad \text{\rm as}\ k\to \infty.
\end{equation*}
Again using assumption \eqref{asm-f}, we find that for fixed $(y,s)\in \R^N\times (s^*,\infty)$,
\begin{equation*}
  \lambda_k^{2p/(p-1)}f\left(\lambda_k^{-2/(p-1)}v_k(y,s)\right)=v_k^p(y,s)s_k^{-p}f(s_k)\to lv^p(y,s)\quad \text{\rm as}\ k\to \infty,
\end{equation*}
where $s_k=\lambda ^{-2/(p-1)}v_k(y,s)$.
Consequently, we derive that $v$ satisfies
\begin{equation*}
  v_t-\Delta v=l v^p+M|\nabla v|^{\frac{2p}{p+1}}\quad \text{\rm in}\ \R^N\times(s^*,\infty).
\end{equation*}
Since the uniqueness of solutions for the corresponding Cauchy problem guarantees $s^*=-\infty$, we deduce that $v$ is a nonnegative  nontrivial solution of
\begin{equation*}\label{lim-v}
  v_t-\Delta v=lv^p+M|\nabla v|^{\frac{2p}{p+1}}\quad \text{\rm in}\ \R^N\times\R.
\end{equation*}
Then under the transformation
\begin{equation*}
  w=l^{\frac{1}{p-1}}v,
\end{equation*}
we know that $w$ is the nonnegative nontrivial solution of the equation
\begin{equation*}
  w_t-\Delta w=w^p+Ml^{-\frac{1}{p+1}}|\nabla w|^{\frac{2p}{p+1}}\quad \text{\rm in}\ \R^N\times \R.
\end{equation*}
When $M\leq \varepsilon_1:=l^{1/(p+1)} \varepsilon_0$, where $\varepsilon_0=\varepsilon_0(N,p)>0$ defined in Theorem \ref{them:int-q=}, the nontrivial property of $w$ contradicts to Theorem \ref{them:int-q=}. This proves Theorem \ref{them:f-g}.
\hfill$\Box$

\vskip 5mm

\noindent{\bf Acknowledgment.} {This work was partially supported by NSFC grants (Nos. 12271423, 12071044), the Shaanxi Fundamental Science Research Project for Mathematics and Physics (No. 23JSY026) and the Fundamental Research Funds for the Central Universities (No. xzy012022005).}

\end{CJK}
\end{document}